\documentclass[12pt,reqno]{article}

\usepackage[usenames]{color}
\usepackage{amssymb}
\usepackage{graphicx}
\usepackage{amscd}

\usepackage[colorlinks=true,
linkcolor=webgreen,
filecolor=webbrown,
citecolor=webgreen]{hyperref}

\definecolor{webgreen}{rgb}{0,.5,0}
\definecolor{recrown}{rgb}{.6,0,0}

\usepackage{color}
\usepackage{fullpage}
\usepackage{float}

\usepackage{graphics,amsmath,amssymb}
\usepackage{amsthm}
\usepackage{amsfonts}
\usepackage{latexsym}
\usepackage{epsf}

\newcommand{\seqnum}[1]{\href{http://www.research.att.com/cgi-bin/access.cgi/as/~njas/sequences/eisA.cgi?Anum=#1}{\underline{#1}}}

\newtheorem{theo}{Theorem}[section]                                                                          %
\newtheorem{lemma}[theo]{Lemma}                                                                              %
\def\dire{dir{\hspace{-0.11em}\raisebox{-0.65ex}{\mbox{\scriptsize $\vec{e}$}}}\,}                           %
\def\stop{\mbox{\footnotesize {\vrule width 6pt height 6pt}}}                                                %
\begin{document}

$\;$

\vspace*{10 mm}                                                 % <-----

% \begin{center}                                                % <-----
% \epsfxsize=4in \leavevmode\epsffile{logo129.eps}              % <----- JIS LOGO
% \end{center}                                                  % <-----

\begin{center}
\vskip 1cm {\LARGE\bf
The Compositions of
the~Differential~Operations and
Gateaux~Directional~Derivative}
\vskip 1cm \large
Branko J. Male\v sevi\' c${\,}^{1)}$\footnote{
${\,}^{1)}\,$This work was supported in part by the project MNTRS, Grant No. ON144020.}
and
Ivana V. Jovovi\' c${\,}^{2)}$\footnote{
${\,}^{2)}\,$PhD student, Faculty of Mathematics, University of Belgrade, Serbia}  \\
University of Belgrade, Faculty of Electrical Engineering                      \\
Bulevar kralja Aleksandra 73, Belgrade, Serbia                                 \\
\href{mailto:malesh@EUnet.yu}{\tt malesh@EUnet.yu}                             \\
\href{mailto:ivana121@EUnet.yu}{\tt ivana121@EUnet.yu}                         \\
\end{center}

\begin{abstract}
\noindent
In this paper we determine the number of the meaningful compositions of higher order
of the differential operations and Gateaux directional derivative.
\end{abstract}

\section{The compositions of the differential operations
of the space $\mathbb R^{\mbox{\footnotesize \textbf\textit{3}}}$}

In the real three-dimensional space $\mathbb R^{3}$ we consider the following sets$:$
\begin{equation}
\mbox{\rm A}_{0}
=
\{ f\!:\! \mathbb R^{3} \!\longrightarrow\! \mathbb R
\, | \, f\!\in\!C^{\infty}(\mathbb R^{3}) \}
\;\;\;
\mbox{and}
\;\;\;
\mbox{\rm A}_{1}
=
\{ \vec{f}\!:\! \mathbb R^{3} \!\longrightarrow\! \mathbb R^{3}
\, | \, \vec{f}\!\in\!\vec{C}^{\infty}(\mathbb R^{3}) \}.
\end{equation}
Then, over the sets $\mbox{\rm A}_{0}$ and $\mbox{\rm A}_{1}$ in the vector analysis,
there are $m=3$ differential operations of the~first-order$:$
\begin{equation}
\begin{array}{l}
\mbox{ \small $ \mbox{\normalsize \rm grad} \, \mbox{\normalsize $f$}
=
\mbox{\normalsize $\nabla_1$} \, \mbox{\normalsize $f$}
\!=\!
\left(\displaystyle \frac{\partial f}{\partial x_1} ,
\displaystyle \frac{\partial f}{\partial x_2} ,
\displaystyle \frac{\partial f}{\partial x_3}\right):
\mbox{\normalsize \rm A}_{0} \longrightarrow \mbox{\normalsize \rm A}_{1} $ },   \\[2.5 ex]
\mbox{ \small $ \mbox{\normalsize \rm curl} \, \vec{\mbox{\normalsize $f$}}
=
\mbox{\normalsize $\nabla_2$} \, \vec{\mbox{\normalsize $f$}}
=
\left( \displaystyle \frac{\partial f_3}{\partial x_2} \!-\!
\displaystyle \frac{\partial f_2}{\partial x_3} ,
       \displaystyle \frac{\partial f_1}{\partial x_3} \!-\!
\displaystyle \frac{\partial f_3}{\partial x_1} ,
       \displaystyle \frac{\partial f_2}{\partial x_1} \!-\!
\displaystyle \frac{\partial f_1}{\partial x_2}   \right) :
\mbox{\normalsize \rm A}_{1} \longrightarrow \mbox{\normalsize \rm A}_{1} $ },   \\[2.5 ex]
\mbox{ \small $ \mbox{\normalsize \rm div} \, \vec{\mbox{\normalsize $f$}}
=
\mbox{\normalsize $\nabla_3$} \, \vec{\mbox{\normalsize $f$}}
=
\displaystyle \frac{\partial f_1}{\partial x_1} \!+\!
\displaystyle \frac{\partial f_2}{\partial x_2} \!+\!
\displaystyle \frac{\partial f_3}{\partial x_3}:
\mbox{\normalsize \rm A}_{1} \longrightarrow \mbox{\normalsize \rm A}_{0} $ }.
\end{array}
\end{equation}
Let us present the number of the meaningful compositions of higher order
over the set ${\cal A}_{3} = \{ \nabla_1, \nabla_2, \nabla_3\}$.
As a well-known fact, there are $m=5$ compositions of the~second-order$:$
\begin{equation}
\begin{array}{l}
\Delta f =  \mbox{div\,grad} \, f = \nabla_3 \circ \nabla_1\, f,   \\[1.5 ex]
\mbox{curl\,curl} \, \vec f = \nabla_2 \circ \nabla_2\, \vec f,    \\[1.5 ex]
\mbox{grad\,div} \, \vec f = \nabla_1 \circ \nabla_3\, \vec f,     \\[1.5 ex]
\mbox{curl\,grad} \, f = \nabla_2 \circ \nabla_1\, f = \vec 0,     \\[1.5 ex]
\mbox{div\,curl} \, \vec f = \nabla_3 \circ \nabla_2\, \vec f = 0.
\end{array}
\end{equation}
Male\v sevi\' c \cite{HiOrd96} proved that there are $m=8$
compositions of the~third-order$:$
\begin{equation}
\begin{array}{l}
\mbox{grad\,div\,grad}\, f = \nabla_1 \circ \nabla_3 \circ \nabla_1\,f,                    \\[1.5 ex]
\mbox{curl\,curl\,curl}\, \vec f = \nabla_2 \circ \nabla_2 \circ \nabla_2\,\vec f,         \\[1.5 ex]
\mbox{div\,grad\,div}\, \vec f = \nabla_3 \circ \nabla_1 \circ \nabla_3\,\vec f,           \\[1.5 ex]
\mbox{curl\,curl\,grad}\, f = \nabla_2 \circ \nabla_2 \circ \nabla_1\,f = \vec 0,          \\[1.5 ex]
\mbox{div\,curl\,grad}\, f = \nabla_3 \circ \nabla_2 \circ \nabla_1\,f = 0,                \\[1.5 ex]
\mbox{div\,curl\,curl}\, \vec f = \nabla_3 \circ \nabla_2 \circ \nabla_2\,\vec f = 0,      \\[1.5 ex]
\mbox{grad\,div\,curl}\, \vec f = \nabla_1 \circ \nabla_3 \circ \nabla_2\,\vec f = \vec 0, \\[1.5 ex]
\mbox{curl\,grad\,div}\, \vec f = \nabla_2 \circ \nabla_1 \circ \nabla_3\,\vec f = \vec 0.
\end{array}
\end{equation}
If we denote by $\mbox{\large \tt f}(k)$ the number of compositions of
the~$k^{\mbox{\scriptsize \rm th}}$-order, then Male\v sevi\' c \cite{HiOrd98} proved$:$
\begin{equation}
\mbox{\large \tt f}(k) = F_{k+3},
\end{equation}
where $F_{k}$ is $k^{\mbox{\scriptsize \rm th}}$ Fibonacci number.

\section{The compositions of the differential operations and
\textsc{Gateaux} directional derivative on the space~$\mathbb R^{\mbox{\footnotesize \textbf\textit{3}}}$}

Let $f \in \mbox{\rm A}_{0}$ be a scalar function and $\vec{e} = (e_1,e_2,e_3)\in \mathbb R^{3}$ be a unit vector.
Thus, the \textit{Gateaux directional derivative} in direction $\vec{e}$ is defined by \cite[p.$\,$71]{Basov05}$:$
\begin{equation}
\mbox{\dire} \, f = \nabla_0 f = \nabla_1 f \cdot \vec{e} =
\frac{\partial f}{\partial x_1} \, e_1 +
\frac{\partial f}{\partial x_2} \, e_2 +
\frac{\partial f}{\partial x_3} \, e_3 :
\mbox{\rm A}_{0} \longrightarrow \mbox{\rm A}_{0}.
\end{equation}
Let us determine the number of the meaningful compositions of higher order over the set
${\cal B}_{3} = \{ \nabla_0, \nabla_1, \nabla_2, \nabla_3 \}$. There exist $m=8$
compositions of the~second-order$:$
\begin{equation}
\label{B_Second}
\begin{array}{l}
           \mbox{\dire\,\dire} \, f = \nabla_0 \circ \nabla_0\, f
= \nabla_1 {\big (} \, \nabla_1 f \cdot \vec{e} \, {\big )} \cdot \vec{e},                  \\[1.5 ex]
           \mbox{grad\,\dire} \, f = \nabla_1 \circ \nabla_0\, f
= \nabla_1 {\big (} \, \nabla_1 f \cdot \vec{e} \, {\big )},                                \\[1.5 ex]
\Delta f = \mbox{div\,grad} \, f = \nabla_3 \circ \nabla_1 \, f,                              \\[1.5 ex]
           \mbox{curl\,curl} \, \vec{f} = \nabla_2 \circ \nabla_2\, \vec{f},                  \\[1.5 ex]
           \mbox{\dire\,div} \, \vec{f} = \nabla_0 \circ \nabla_3\, \vec{f}
= {\big (}\nabla_1 \circ \nabla_3 \vec{f}{\big )} \cdot \vec{e},                          \\[1.5 ex]
           \mbox{grad\,div} \, \vec{f} = \nabla_1 \circ \nabla_3\, \vec{f},                   \\[1.5 ex]
\mbox{curl\,grad} \, f = \nabla_2 \circ \nabla_1\, f = \vec{0},                               \\[1.5 ex]
           \mbox{div\,curl} \, \vec{f} = \nabla_3 \circ \nabla_2\, \vec{f} = 0;
\end{array}
\end{equation}
that is, there exist $m=16$ compositions of the~third-order$:$
\begin{equation}
\begin{array}{l}
\mbox{\dire\,\dire\,\dire}\, f = \nabla_0 \circ \nabla_0 \circ \nabla_0\, f,                   \\[1.5 ex]
\mbox{grad\,\dire\,\dire}\, f = \nabla_1 \circ \nabla_0 \circ \nabla_0\, f,                    \\[1.5 ex]
\mbox{div\,grad\,\dire}\, f = \nabla_3 \circ \nabla_1 \circ \nabla_0\, f,                      \\[1.5 ex]
\mbox{\dire\,div\,grad}\, f = \nabla_0 \circ \nabla_3 \circ \nabla_1\, f,                      \\[1.5 ex]
\mbox{grad\,div\,grad}\, f = \nabla_1 \circ \nabla_3 \circ \nabla_1\, f,                       \\[1.5 ex]
\mbox{curl\,curl\,curl}\, \vec{f} = \nabla_2 \circ \nabla_2 \circ \nabla_2\, \vec{f},          \\[1.5 ex]
\mbox{\dire\,\dire\,div}\, \vec{f} = \nabla_0 \circ \nabla_0 \circ \nabla_3\, \vec{f},         \\[1.5 ex]
\mbox{grad\,\dire\,div}\, \vec{f} = \nabla_1 \circ \nabla_0 \circ \nabla_3\, \vec{f},          \\[1.5 ex]
\mbox{div\,grad\,div}\, \vec{f} = \nabla_3 \circ \nabla_1 \circ \nabla_3\, \vec{f},            \\[1.5 ex]
\mbox{curl\,grad\,\dire}\, f = \nabla_2 \circ \nabla_1 \circ\nabla_0\, \vec{f} = \vec{0},      \\[1.5 ex]
\mbox{curl\,curl\,grad}\, f = \nabla_2 \circ \nabla_2 \circ \nabla_1\, f = \vec{0},            \\[1.5 ex]
\mbox{div\,curl\,grad}\, f = \nabla_3 \circ \nabla_2 \circ \nabla_1\, f = 0,                   \\[1.5 ex]
\mbox{div\,curl\,curl}\, \vec{f} = \nabla_3 \circ \nabla_2 \circ \nabla_2\, \vec{f} = 0,       \\[1.5 ex]
\mbox{\dire\,div\,curl}\, \vec{f} = \nabla_0 \circ \nabla_3 \circ \nabla_2\, \vec{f} = 0,      \\[1.5 ex]
\mbox{grad\,div\,curl}\, \vec{f} = \nabla_1 \circ \nabla_3 \circ \nabla_2\, \vec{f} = \vec{0}, \\[1.5 ex]
\mbox{curl\,grad\,div}\, \vec{f} = \nabla_2 \circ \nabla_1 \circ \nabla_3\, \vec{f} = \vec{0}.
\end{array}
\end{equation}

\break

\noindent
Using the method from the paper \cite{HiOrd98} let us define a binary relation $\sigma$
``{\em to be in composition}''$:$ \mbox{$\nabla_{i} \,\sigma\, \nabla_{j} = \top$}
iff the composition $\nabla_{j} \circ \nabla_{i}$ is meaningful.
Thus, Cayley table of the relation $\sigma$ is determined with
\begin{equation}
\label{Cayley_4}
\begin{array}{c|cccc}
\sigma & \nabla_{0} & \nabla_{1} & \nabla_{2} & \nabla_{3}  \\ \hline
\nabla_{0} & \top & \top       & \bot       & \bot        \\
\nabla_{1} & \bot & \bot       & \top       & \top        \\
\nabla_{2} & \bot & \bot       & \top       & \top        \\
\nabla_{3} & \top & \top       & \bot       & \bot
\end{array}
\end{equation}

\medskip
\noindent
Let us form the graph according to the following rule$:$
if $\nabla_{i} \,\sigma\, \nabla_{j} = \top$ let vertex $\nabla_{j}$ be under vertex $\nabla_{i}$
and let there exist an edge from the vertex $\nabla_{i}$ to the vertex $\nabla_{j}$. Further on,
let us denote by $\nabla_{\!\!-1}$ nowhere-defined function $\vartheta$, where domain and range are
the empty sets \cite{HiOrd96}. We shall define $\nabla_{\!\!-1} \,\sigma\, \nabla_{i} = \top$ $(i\!=\!0,1,2,3,4)$.
For the set ${\cal B}_{3} \cup \{ \nabla_{\!\!-1} \}$ the graph of the walks, determined previously,
is a tree with the root in the vertex $\nabla_{\!\!-1}$.

%*********************************************************************** Fig 1 **
\hspace*{21.0 mm}                                                              %*
\setlength{\unitlength}{0.17 cc}                                               %*
\begin{picture}(150,50)(0,0)                                                   %*
\thicklines                                                                    %*
\put(55,40){\circle*{0.8}}                                                     %*
\put(56,41){\scriptsize$\nabla_{\!\!-1}$}                                      %*
\put(130,40){\scriptsize $\mbox{\footnotesize \tt g}(0)=\;1$}                  %*
                                                                               %*
\put(5,30){\line(5,1){50}}                                                     %*
\put(5,30){\circle*{0.8}}                                                      %*
\put(3,32){\scriptsize$\nabla_{0}$}                                            %*
\put(35,30){\line(2,1){20}}                                                    %*
\put(35,30){\circle*{0.8}}                                                     %*
\put(32,32){\scriptsize$\nabla_{1}$}                                           %*
\put(75,30){\line(-2,1){20}}                                                   %*
\put(75,30){\circle*{0.8}}                                                     %*
\put(75,32){\scriptsize$\nabla_{2}$}                                           %*
\put(105,30){\line(-5,1){50}}                                                  %*
\put(105,30){\circle*{0.8}}                                                    %*
\put(105,32){\scriptsize$\nabla_{3}$}                                          %*
\put(130,30){\scriptsize $\mbox{\footnotesize \tt g}(1)=\;4$}                  %*
\put(-5,20){\line(1,1){10}}                                                    %*
\put(-5,20){\circle*{0.8}}                                                     %*
\put(-8,22){\scriptsize$\nabla_{0}$}                                           %*
\put(15,20){\line(-1,1){10}}                                                   %*
\put(15,20){\circle*{0.8}}                                                     %*
\put(15,22){\scriptsize$\nabla_{1}$}                                           %*
\thinlines                                                                     %*
\put(25,20){\line(1,1){10}}                                                    %*
\put(25,20){\circle*{0.8}}                                                     %*
\put(22,22){\scriptsize$\nabla_{2}$}                                           %*
\thicklines                                                                    %*
\put(45,20){\line(-1,1){10}}                                                   %*
\put(45,20){\circle*{0.8}}                                                     %*
\put(43,22){\scriptsize$\nabla_{3}$}                                           %*
\put(65,20){\line(1,1){10}}                                                    %*
\put(65,20){\circle*{0.8}}                                                     %*
\put(62,22){\scriptsize$\nabla_{2}$}                                           %*
\thinlines                                                                     %*
\put(85,20){\line(-1,1){10}}                                                   %*
\put(85,20){\circle*{0.8}}                                                     %*
\put(85,22){\scriptsize$\nabla_{3}$}                                           %*
\thicklines                                                                    %*
\put(95,20){\line(1,1){10}}                                                    %*
\put(95,20){\circle*{0.8}}                                                     %*
\put(92,22){\scriptsize$\nabla_{0}$}                                           %*
\put(115,20){\line(-1,1){10}}                                                  %*
\put(115,20){\circle*{0.8}}                                                    %*
\put(115,22){\scriptsize$\nabla_{1}$}                                          %*
\put(130,20){\scriptsize $\mbox{\footnotesize \tt g}(2)=\;8$}                  %*
                                                                               %*
\put(-7.5,15){\line(1,2){2.5}}                                                 %*
\put(-2.5,15){\line(-1,2){2.5}}                                                %*
\thinlines                                                                     %*
\put(12.5,15){\line(1,2){2.5}}                                                 %*
\thicklines                                                                    %*
\put(17.5,15){\line(-1,2){2.5}}                                                %*
\thinlines                                                                     %*
\put(22.5,15){\line(1,2){2.5}}                                                 %*
\put(27.5,15){\line(-1,2){2.5}}                                                %*
\thicklines                                                                    %*
\put(42.5,15){\line(1,2){2.5}}                                                 %*
\put(47.5,15){\line(-1,2){2.5}}                                                %*
\put(62.5,15){\line(1,2){2.5}}                                                 %*
\thinlines                                                                     %*
\put(67.5,15){\line(-1,2){2.5}}                                                %*
\put(82.5,15){\line(1,2){2.5}}                                                 %*
\put(87.5,15){\line(-1,2){2.5}}                                                %*
\thicklines                                                                    %*
\put(92.5,15){\line(1,2){2.5}}                                                 %*
\put(97.5,15){\line(-1,2){2.5}}                                                %*
\thinlines                                                                     %*
\put(112.5,15){\line(1,2){2.5}}                                                %*
\thicklines                                                                    %*
\put(117.5,15){\line(-1,2){2.5}}                                               %*
\put(130,15){\scriptsize $\mbox{\footnotesize \tt g}(3)=\;16$}                 %*
\thinlines                                                                     %*
\put(55,4.5){\small Fig. $1$}                                                  %*
\end{picture} %******************************************************************

\noindent
Let $\mbox{\large \tt g}(k)$ be the number of the meaningful compositions
of the $k^{\mbox{\scriptsize \rm th}}$-order of the functions from ${\cal B}_{3}$.
Let $\mbox{\large \tt g}_{i}(k)$ be the number of the meaningful compositions of the
$k^{\mbox{\scriptsize \rm th}}$-order beginning from the left by $\nabla_{i}$.
Then $\mbox{\large \tt g}(k)
=
\mbox{\large \tt g}_{0}(k)
+
\mbox{\large \tt g}_{1}(k)
+
\mbox{\large \tt g}_{2}(k)
+
\mbox{\large \tt g}_{3}(k)$.
Based on the partial self similarity of the tree (Fig.\,$1$) we get equalities
\begin{equation}
\begin{array}{l}
\mbox{\large \tt g}_{0}(k)
=
\mbox{\large \tt g}_{0}(k-1)
+
\mbox{\large \tt g}_{1}(k-1),                                                  \\[1.5 ex]
\mbox{\large \tt g}_{1}(k)
=
\mbox{\large \tt g}_{2}(k-1)
+
\mbox{\large \tt g}_{3}(k-1),                                                  \\[1.5 ex]
\mbox{\large \tt g}_{2}(k)
=
\mbox{\large \tt g}_{2}(k-1)
+
\mbox{\large \tt g}_{3}(k-1),                                                  \\[1.5 ex]
\mbox{\large \tt g}_{3}(k)
=
\mbox{\large \tt g}_{0}(k-1)
+
\mbox{\large \tt g}_{1}(k-1).
\end{array}
\end{equation}
Hence, a recurrence for $\mbox{\large \tt g}(k)$ can be derived as follows$:$
\begin{equation}
\mbox{\large \tt g}(k) \!=\! 2 \, \mbox{\large \tt g}(k-1).
\end{equation}
Based on the initial value $\mbox{\large \tt g}(1)=4$,
we can conclude$:$
\begin{equation}
\mbox{\large \tt g}(k) = 2^{k+1}.
\end{equation}

\section{The compositions of the differential operations
of the space $\mathbb R^{\mbox{\footnotesize \textbf\textit{n}}}$}

Let us present the number of the meaningful compositions of differential operations in the vector analysis
of the space $\mathbb R^{n}$, where differential operations $\nabla_{r}$ $(r \!=\! 1,\ldots,n)$ are
defined over non-empty corresponding sets $\mbox{A}_{s}$ $(s \!=\! 1,\ldots,m$ and
$m \!=\! \lfloor n/2 \rfloor$, $n \!\geq\! 3)$
according to the papers \cite{HiOrd98}, \cite{HiOrd06}$:$

\vspace*{-5.0 mm}

\begin{equation}
\label{A_12}{}
\begin{tabular}{cc}
$\begin{array}{ll}
\mbox{\small $\mbox{$\cal A$}_{n}\;(n\!=\!2m)$:}\!\!
           & \mbox{\small $\nabla_{1}$}   : \mbox{A}_{0} \!\rightarrow\! \mbox{A}_{1}   \\
           & \mbox{\small $\nabla_{2}$}   : \mbox{A}_{1} \!\rightarrow\! \mbox{A}_{2}   \\
           & \,\,\vdots                                                                 \\
           & \mbox{\small $\nabla_{i}$}   : \mbox{A}_{i-1} \!\rightarrow\! \mbox{A}_{i} \\
           & \,\,\vdots                                                                 \\
           & \mbox{\small $\nabla_{m}$}   : \mbox{A}_{m-1} \!\rightarrow\! \mbox{A}_{m} \\
           & \mbox{\small $\nabla_{m+1}$} : \mbox{A}_{m} \!\rightarrow\! \mbox{A}_{m-1} \\
           & \,\,\vdots                                                                 \\
           & \mbox{\small $\nabla_{n-j}$} : \mbox{A}_{j+1} \!\rightarrow\! \mbox{A}_{j} \\
           & \,\,\vdots                                                                 \\
           & \mbox{\small $\nabla_{n-1}$} : \mbox{A}_{2} \!\rightarrow\! \mbox{A}_{1}   \\
           & \mbox{\small $\nabla_{n}$}   : \mbox{A}_{1} \!\rightarrow\! \mbox{A}_{0}
             \mbox{\normalsize ,}
\end{array}$
           &
$ \begin{array}{ll}
\mbox{\small $\mbox{$\cal A$}_{n}\;(n\!=\!2m\!+\!1)$:}\!\!
           & \mbox{\small $\nabla_{1}$}   : \mbox{A}_{0} \!\rightarrow\! \mbox{A}_{1}   \\
           & \mbox{\small $\nabla_{2}$}   : \mbox{A}_{1} \!\rightarrow\! \mbox{A}_{2}   \\
           & \,\,\vdots                                                                 \\
           & \mbox{\small $\nabla_{i}$}   : \mbox{A}_{i-1} \!\rightarrow\! \mbox{A}_{i} \\
           & \,\,\vdots                                                                 \\
           & \mbox{\small $\nabla_{m}$}   : \mbox{A}_{m-1} \!\rightarrow\! \mbox{A}_{m} \\
           & \mbox{\small $\nabla_{m+1}$} : \mbox{A}_{m} \!\rightarrow\! \mbox{A}_{m}   \\
           & \mbox{\small $\nabla_{m+2}$} : \mbox{A}_{m} \!\rightarrow\! \mbox{A}_{m-1} \\
           & \,\,\vdots                                                                 \\
           & \mbox{\small $\nabla_{n-j}$} : \mbox{A}_{j+1} \!\rightarrow\! \mbox{A}_{j} \\
           & \,\,\vdots                                                                 \\
           & \mbox{\small $\nabla_{n-1}$} : \mbox{A}_{2} \!\rightarrow\! \mbox{A}_{1}   \\
           & \mbox{\small $\nabla_{n}$}   : \mbox{A}_{1} \!\rightarrow\! \mbox{A}_{0}
             \mbox{\normalsize .}
\end{array}$
\end{tabular}
\end{equation}

\vspace*{-1.0 mm}

\noindent
Let us define {\em higher order differential operations} as the meaningful compositions of higher
order of differential operations from the set ${\cal A}_{n} = \{\nabla_{1}, \dots, \nabla_{n}\}$.
The number of the higher order differential operations is given according to the paper \cite{HiOrd98}.
Let us define a binary relation $\rho$ ``{\em to be in composition}''$:$
$\nabla_{i} \,\rho\, \nabla_{j} = \top$ iff the composition $\nabla_{j} \circ \nabla_{i}$
is meaningful. Thus, Cayley table of the relation $\rho$ is determined with
\begin{equation}
\label{Rho}
\mbox{\normalsize $\nabla_{i} \,\rho\, \nabla_{j}$}
=
\left\{
\begin{array}{lll}
\top &,& (j = i + 1)    \vee   (i + j = n + 1);                                 \\[1.0 ex]
\bot &,& \mbox{\normalsize otherwise}.
\end{array}
\right.
\end{equation}
Let us form the adjacency matrix $\mbox{\large \tt A} = [a_{ij}]\in\{\,0,1\}^{n\times n}$
associated with the graph, which is determined by the relation $\rho$. Thus, according
to the paper \cite{HiOrd06}, the following statement is true.
\begin{theo}
\label{Th_3_1}
Let $P_{n}(\lambda) \!=\! |\mbox{\large \tt A} \!-\! \lambda \mbox{\large \tt I}| \!=\!
\alpha_{0} \lambda^{n} + \alpha_{1} \lambda^{n-1} + \dots + \alpha_{n}$
be the characteristic polynomial of the matrix $\mbox{\large \tt A}$
and  $v_{n} = [ \, 1 \, \dots \, 1 \, ]_{1 \times n}$.
If we denote by $\mbox{\large \tt f}(k)$ the number of the
$k^{\it \footnotesize th}\!$-order differential operations,
then the following formulas are true$:$
\begin{equation}
\mbox{\large \tt f}(k) = v_n \cdot \mbox{\large \tt A}^{k-1} \cdot v^{T}_n
\end{equation}
and
\begin{equation}
\alpha_{0} \mbox{\large \tt f}(k) + \alpha_{1} \mbox{\large \tt f}(k-1) + \dots + \alpha_{n} \mbox{\large \tt f}(k-n) = 0
\quad (k > n).
\end{equation}
\end{theo}

\break

\begin{lemma}
\label{lemma_stara_1}
Let $P_{n}(\lambda)$ be the characteristic polynomial of the matrix $\mbox{\large \tt A}$.
Then the following recurrence is true$:$
\begin{equation}
\label{lemma_stara_1_Form_1}
P_{n}(\lambda) = \lambda^2 {\big (} P_{n-2}(\lambda) - P_{n-4}(\lambda) {\big )}.
\end{equation}
\end{lemma}
\begin{lemma}
\label{lemma_stara_2}
Let $P_{n}(\lambda)$ be the characteristic polynomial of the matrix $\mbox{\large \tt A}$.
Then it has the following explicit representation:
\begin{equation}
\label{lemma_stara_2_Form_2}
\quad
P_{n}(\lambda)
=
\left\{
\begin{array}{ccl}
\displaystyle\sum\limits_{k=1}^{\lfloor \frac{n+2}{4} \rfloor +1}{(-1)^{k-1}
{ \: \mbox{\scriptsize $\displaystyle\frac{n}{2}\!-\!k\!+\!2$} \:
\choose
\: \mbox{\scriptsize $k\!-\!1$} \: }
\lambda^{n-2k+2}}
\!\!&\!\!,\!\!& n\!=\!2m;                                              \\[2.0 ex]
\!\!\!\displaystyle\sum\limits_{k=1}^{\lfloor \frac{n+2}{4} \rfloor +2}{\!\!\!\!(-1)^{k-1}\!{\Bigg (}
\!{ \: \mbox{\scriptsize $\displaystyle\frac{n\!+\!3}{2}\!-\!k$} \:
\choose
\: \mbox{\scriptsize $k\!-\!1$} \: }
\!+\!
{ \: \mbox{\scriptsize $\displaystyle\frac{n\!+\!3}{2}\!-\!k$} \:
\choose
\: \mbox{\scriptsize $k\!-\!2$} \: } \! \lambda \! {\Bigg )} \lambda^{n-2k+2}}
\!&\!\!,\!\!& n\!=\!2m\!+\!1.\!\!\!\!
\end{array}
\right.
\end{equation}
\end{lemma}

\noindent
The number of the higher order differential operations is determined by
corresponding recurrence, which for dimension $n = 3, 4, 5,\dots, 10\,$,
we refer according to \cite{HiOrd98}:

{ \footnotesize
\begin{center}
\begin{tabular}{|c|c|}
\hline {\small \rm Dimension}:
 &     {\small \rm \quad Recurrence for the number of the $k^{\mbox{\scriptsize \rm th}}$-order differential operations: \quad\quad} \\ \hline
$n = \;$ 3 & $ \mbox{\normalsize \tt f}(k)
             = \mbox{\normalsize \tt f}(k-1)
             + \mbox{\normalsize \tt f}(k-2)$                                         \\ \hline
$n = \;$ 4 & $\mbox{\normalsize \tt f}(k)
             = 2 \mbox{\normalsize \tt f}(k-2)$                                       \\ \hline
$n = \;$ 5 & $\mbox{\normalsize \tt f}(k)
             = \mbox{\normalsize \tt f}(k-1)
             + 2 \mbox{\normalsize \tt f}(k-2)
             - \mbox{\normalsize \tt f}(k-3)$                                         \\ \hline
$n = \;$ 6 & $\mbox{\normalsize \tt f}(k)
             = 3 \mbox{\normalsize \tt f}(k-2)
             - \mbox{\normalsize \tt f}(k-4)$                                         \\ \hline
$n = \;$ 7 & $\mbox{\normalsize \tt f}(k)
             = \mbox{\normalsize \tt f}(k-1)
             + 3 \mbox{\normalsize \tt f}(k-2)
             - 2 \mbox{\normalsize \tt f}(k-3)
             - \mbox{\normalsize \tt f}(k-4)$                                       \\ \hline
$n = \;$ 8 & $\mbox{\normalsize \tt f}(k)
             = 4 \mbox{\normalsize \tt f}(k-2)
             - 3 \mbox{\normalsize \tt f}(k-4)$                                     \\ \hline
$n = \;$ 9 & $\mbox{\normalsize \tt f}(k)
             = \mbox{\normalsize \tt f}(k-1)
             + 4 \mbox{\normalsize \tt f}(k-2)
             - 3 \mbox{\normalsize \tt f}(k-3)
             - 3 \mbox{\normalsize \tt f}(k-4)
             + \mbox{\normalsize \tt f}(k-5)$                                      \\ \hline
$n =   $ 10& $\mbox{\normalsize \tt f}(k)
             = 5 \mbox{\normalsize \tt f}(k-2)
             - 6 \mbox{\normalsize \tt f}(k-4)
             + \mbox{\normalsize \tt f}(k-6)$                                      \\ \hline
\end{tabular}
\end{center} }

\smallskip
\noindent
For considered dimensions $n=3,4,5, \dots ,10$, the values of the function $\mbox{\large \tt f}(k)$,
for small values of the argument $k$, are given in the database of integer sequences \cite{Sloane07}
as sequences
\seqnum{A020701} $(n=3)$, \seqnum{A090989} $(n=4)$,
\seqnum{A090990} $(n=5)$, \seqnum{A090991} $(n=6)$,
\seqnum{A090992} $(n=7)$, \seqnum{A090993} $(n=8)$,
\seqnum{A090994} $(n=9)$, \seqnum{A090995} $(n=10)$, respectively.

\section{The compositions of the differential operations
and Gateaux directional derivative of the space $\mathbb R^{\mbox{\footnotesize \textbf\textit{n}}}$}

Let $f \in A_{0}$ be a scalar function and $\vec{e} = (e_1,\dots,e_n) \in \mathbb R^{n}$
be a unit vector. Thus, the \textit{Gateaux directional derivative} in direction $\vec{e}$
is defined by \cite[p.$\,$71]{Basov05}$:$
\begin{equation}
\mbox{\dire} \, f = \nabla_0 f =
\displaystyle \sum\limits_{k=1}^{n}{ \frac{\partial f}{\partial x_k} \, e_k} : A_{0} \longrightarrow A_{0}.
\end{equation}
Let us extend the set of differential operations ${\cal A}_{n} = \{ \nabla_{1},\dots,\nabla_{n} \}$
with Gateaux directional derivational to the set
${\cal B}_{n} = {\cal A}_{n} \cup \{ \nabla_{0} \} = \{ \nabla_{0},\nabla_{1},\dots,\nabla_{n} \}$$:$
\begin{equation}
\label{B_12}
\mbox{$\begin{array}{ll}
\mbox{\small $\mbox{$\cal B$}_{n}\;(n\!=\!2m)$:}\!\!
             & \mbox{\small $\nabla_{0}$}   : \mbox{A}_{0} \!\rightarrow\! \mbox{A}_{0}   \\
             & \mbox{\small $\nabla_{1}$}   : \mbox{A}_{0} \!\rightarrow\! \mbox{A}_{1}   \\
             & \mbox{\small $\nabla_{2}$}   : \mbox{A}_{1} \!\rightarrow\! \mbox{A}_{2}   \\
             & \,\,\vdots                                                                 \\
             & \mbox{\small $\nabla_{i}$}   : \mbox{A}_{i-1} \!\rightarrow\! \mbox{A}_{i} \\
             & \,\,\vdots                                                                 \\
             & \mbox{\small $\nabla_{m}$}   : \mbox{A}_{m-1} \!\rightarrow\! \mbox{A}_{m} \\
             & \mbox{\small $\nabla_{m+1}$} : \mbox{A}_{m} \!\rightarrow\! \mbox{A}_{m-1} \\
             & \,\,\vdots                                                                 \\
             & \mbox{\small $\nabla_{n-j}$} : \mbox{A}_{j+1} \!\rightarrow\! \mbox{A}_{j} \\
             & \,\,\vdots                                                                 \\
             & \mbox{\small $\nabla_{n-1}$} : \mbox{A}_{2} \!\rightarrow\! \mbox{A}_{1}   \\
             & \mbox{\small $\nabla_{n}$}   : \mbox{A}_{1} \!\rightarrow\! \mbox{A}_{0}
             \mbox{\normalsize ,}
\end{array} $}
\quad
\mbox{$ \begin{array}{ll}
\mbox{\small $\mbox{$\cal B$}_{n}\;(n\!=\!2m\!+\!1)$:}\!\!
             & \mbox{\small $\nabla_{0}$}   : \mbox{A}_{0} \!\rightarrow\! \mbox{A}_{0}   \\
             & \mbox{\small $\nabla_{1}$}   : \mbox{A}_{0} \!\rightarrow\! \mbox{A}_{1}   \\
             & \mbox{\small $\nabla_{2}$}   : \mbox{A}_{1} \!\rightarrow\! \mbox{A}_{2}   \\
             & \,\,\vdots                                                                 \\
             & \mbox{\small $\nabla_{i}$}   : \mbox{A}_{i-1} \!\rightarrow\! \mbox{A}_{i} \\
             & \,\,\vdots                                                                 \\
             & \mbox{\small $\nabla_{m}$}   : \mbox{A}_{m-1} \!\rightarrow\! \mbox{A}_{m} \\
             & \mbox{\small $\nabla_{m+1}$} : \mbox{A}_{m} \!\rightarrow\! \mbox{A}_{m}   \\
             & \mbox{\small $\nabla_{m+2}$} : \mbox{A}_{m} \!\rightarrow\! \mbox{A}_{m-1} \\
             & \,\,\vdots                                                                 \\
             & \mbox{\small $\nabla_{n-j}$} : \mbox{A}_{j+1} \!\rightarrow\! \mbox{A}_{j} \\
             & \,\,\vdots                                                                 \\
             & \mbox{\small $\nabla_{n-1}$} : \mbox{A}_{2} \!\rightarrow\! \mbox{A}_{1}   \\
             & \mbox{\small $\nabla_{n}$}   : \mbox{A}_{1} \!\rightarrow\! \mbox{A}_{0}
             \mbox{\normalsize .}
\end{array} $}
\end{equation}
Let us define \textit{higher order differential operations with Gateaux derivative}
as the meaningful compositions of higher order of the  functions from the set
${\cal B}_{n} = \{ \nabla_{0},\nabla_{1},\dots, \nabla_{n} \}$.
We determine the number of the higher order differential operations with Gateaux derivative
by defining a binary relation~$\sigma$ ``{\em to~be~in~composition}''$:$
\begin{equation}
\label{Rho'2}
\nabla_{i} \,\sigma\, \nabla_{j}
=
\left\{
\begin{array}{lll}
\top\! &\!\!,\!\!& (i\!=\!0 \wedge j\!=\!0) \vee (i\!=\!n \wedge j\!=\!0)
\vee (j\!=\!i\!+\!1) \vee (i\!+\!j\!=\!n\!+\!1);               \\[1.0 ex]
\bot\! &\!\!,\!\!& \mbox{\normalsize otherwise}.
\end{array}
\right.
\end{equation}
Let us form the adjacency matrix $\mbox{\large \tt B} = [b_{ij}]\in\{\,0,1\}^{(n+1) \times n}$
associated with the graph, which is determined by relation $\sigma$. Thus, analogously
to the paper \cite{HiOrd06}, the following statement is true.
\begin{theo}
\label{Th_4_1}
Let $Q_{n}(\lambda) \!=\! |\mbox{\large \tt B} \!- \! \lambda \mbox{\large \tt I}| =
\beta_{0} \lambda^{n+1} + \beta_{1}\lambda^{n} + \dots + \beta_{n+1}$
be the characteristic polynomial of the matrix $\mbox{\large \tt B}$
and  $v_{n+1} = [ \, 1 \, \dots \, 1 \, ]_{1  \times (n+1)}$.
If we denote by $\mbox{\large \tt g}(k)$ the number of the
$k^{\it \footnotesize th}\!$-order differential operations with Gateaux derivative,
then the following formulas are true$:$
\begin{equation}
\mbox{\large \tt g}(k) = v_{n+1} \cdot \mbox{\large \tt B}^{k-1} \cdot v^{T}_{n+1}
\end{equation}
and
\begin{equation}
\beta_{0} \mbox{\large \tt g}(k)
+
\beta_{1} \mbox{\large \tt g}(k-1)
+
\dots
+
\beta_{n+1} \mbox{\large \tt g}(k-(n+1)) = 0
\quad (k > n\!+\!1).
\end{equation}
\end{theo}
\begin{lemma}
\label{lemma_4_2}
Let $Q_{n}(\lambda)$ and $P_{n}(\lambda)$ be the characteristic polynomials
of the matrices $\mbox{\large \tt B}$ and $\mbox{\large \tt A}$ respectively.
Then the following equality is true$:$
\begin{equation}
\label{veza izmedju A i B}
Q_{n}(\lambda) = \lambda^2 P_{n-2}(\lambda) - \lambda P_{n}(\lambda).
\end{equation}
\end{lemma}

\noindent
{\bf Proof.}
Let us determine the characteristic polynomial $Q_{n}(\lambda) = |\mbox{\large \tt B} - \lambda \mbox{\large \tt I}|$ by

\begin{equation}
Q_{n}(\lambda) =
\mbox{\footnotesize $\left|
\begin{array}{rrrrrrrrr}
1-\lambda & 1 & 0 & 0 & \dots & 0 & 0 & 0 & 0 \\
0 & -\lambda  & 1 & 0 & \dots & 0 & 0 & 0 & 1 \\
0 & 0 & -\lambda  & 1 & \dots & 0 & 0 & 1 & 0 \\
\vdots & \vdots & \vdots & \vdots & \ddots &
\vdots & \vdots & \vdots & \vdots             \\
0 & 0 & 0 & 1 & \dots & 0 & -\lambda & 1 & 0  \\
0 & 0 & 1 & 0 & \dots & 0 & 0 & -\lambda & 1  \\
1 & 1 & 0 & 0 & \dots & 0 & 0 & 0 & -\lambda
\end{array}
\right|$}\,.
\end{equation}

\noindent
Expanding the determinant $Q_{n}(\lambda)$
by the first column we have
\begin{equation}
\label{jed. 1}
Q_{n}(\lambda)
=
(1-\lambda) P_{n}(\lambda) + (-1)^{n+2} D_{n}(\lambda),
\end{equation}
where is

\vspace*{-1.0 mm}

\begin{equation}
\label{detD 1}
D_{n}(\lambda) =
\mbox{\footnotesize $\left|
\begin{array}{rrrrrrrrr}
1        & 0        & 0 & 0 & \dots & 0 & 0 & 0 & 0 \\
-\lambda & 1        & 0 & 0 & \dots & 0 & 0 & 0 & 1 \\
0        & -\lambda & 1 & 0 & \dots & 0 & 0 & 1 & 0 \\
\vdots & \vdots & \vdots & \vdots & \ddots &
\vdots & \vdots & \vdots & \vdots            \\
0 & 0 & 0 & 1 & \dots & -\lambda & 1 & 0 & 0 \\
0 & 0 & 1 & 0 & \dots & 0 & -\lambda & 1 & 0 \\
0 & 1 & 0 & 0 & \dots & 0 & 0 & -\lambda & 1
\end{array}
\right|$}\,.
\end{equation}

\noindent
Let us expand the determinant $D_{n}(\lambda)$ by the first row and then,
in the next step, let us multiply the first row by $-1$ and add it to the last row.
Then, we obtain the determinant of order $n-1:$

\begin{equation}
\label{detD 2}
D_{n}(\lambda) =
\mbox{\footnotesize $\left|
\begin{array}{rrrrrrrrr}
       1 & 0 & 0 & 0 & \dots & 0 & 0 & 0 & 1 \\
-\lambda & 1 & 0 & 0 & \dots & 0 & 0 & 1 & 0 \\
0 & -\lambda & 1 & 0 & \dots & 0 & 1 & 0 & 0 \\
\vdots & \vdots & \vdots & \vdots & \ddots &
\vdots & \vdots & \vdots & \vdots            \\
0 & 0 & 1 & 0 & \dots & -\lambda & 1 & 0 & 0 \\
0 & 1 & 0 & 0 & \dots & 0 & -\lambda & 1 & 0 \\
0 & 0 & 0 & 0 & \dots & 0 & 0 & -\lambda & 0
\end{array}
\right|$}\,.
\end{equation}

\noindent
Expanding the previous determinant by the last column we have

\begin{equation}
\label{detD 3}
D_{n}(\lambda) =
(-1)^{n} \mbox{\footnotesize $\left|
\begin{array}{rrrrrrrrr}
-\lambda & 1 & 0 & 0 & \dots & 0 & 0 & 0 & 1  \\
0 & -\lambda & 1 & 0 & \dots & 0 & 0 & 1 & 0  \\
0 & 0 & -\lambda & 1 & \dots & 0 & 1 & 0 & 0  \\
\vdots & \vdots & \vdots & \vdots & \ddots &
\vdots & \vdots & \vdots & \vdots             \\
0 & 0 & 1 & 0 & \dots & 0 & -\lambda & 1 & 0  \\
0 & 1 & 0 & 0 & \dots & 0 & 0 & -\lambda & 1  \\
0 & 0 & 0 & 0 & \dots & 0 & 0 & 0 & -\lambda
\end{array}
\right|$}\,.
\end{equation}

\break

\noindent
If we expand the previous determinant by the last row, and if we expand
the obtained determinant by the first column, we have the determinant of order $n-4:$

\begin{equation}
\label{detD 4}
D_{n}(\lambda) =
(-1)^{n}\lambda^2 \mbox{\footnotesize $\left|
\begin{array}{rrrrrrrrr}
-\lambda & 1 & 0 & 0 & \dots & 0 & 0 & 0 & 1  \\
0 & -\lambda & 1 & 0 & \dots & 0 & 0 & 1 & 0  \\
0 & 0 & -\lambda & 1 & \dots & 0 & 1 & 0 & 0  \\
\vdots & \vdots & \vdots & \vdots & \ddots &
\vdots & \vdots &\vdots & \vdots              \\
0 & 0 & 1 & 0 & \dots & 0 & -\lambda & 1 & 0  \\
0 & 1 & 0 & 0 & \dots & 0 & 0 & -\lambda & 1  \\
1 & 0 & 0 & 0 & \dots & 0 & 0 & 0 & -\lambda
\end{array}
\right|$}\,.
\end{equation}

\noindent
In other words

\begin{equation}
\label{jedD 1}
D_{n}(\lambda) =(-1)^{n}\lambda^2P_{n-4}(\lambda) .
\end{equation}

\noindent
>From equalities {\rm (\ref{jedD 1})} and {\rm (\ref{jed. 1})} there follows$:$
\begin{equation}
Q_{n}(\lambda) =(1-\lambda) P_{n}(\lambda) + \lambda^2P_{n-4}(\lambda) .
\end{equation}

\noindent
On the basis of Lemma {\rm \ref{lemma_stara_1}.} the following equality is true$:$

\begin{equation}
Q_{n}(\lambda) = \lambda^2 P_{n-2}(\lambda) -\lambda P_{n}(\lambda) .\;~\stop
\end{equation}

\bigskip

\begin{lemma}
\label{lemma_4_3}
Let $Q_{n}(\lambda)$ be the characteristic polynomial of the matrix $\mbox{\large \tt B}$.
Then the following recurrence is true$:$
\begin{equation}
Q_{n}(\lambda) = \lambda^2 {\big (} Q_{n-2}(\lambda) - Q_{n-4}(\lambda) {\big )}.
\end{equation}
\end{lemma}

\noindent
{\bf Proof.}
On the basis of Lemma {\rm \ref{lemma_stara_1}.} and Lemma {\rm \ref{lemma_4_2}.}
there follows the statement.~\stop

\bigskip

\begin{lemma}
\label{lemma_4_4}
Let $Q_{n}(\lambda)$ be the characteristic polynomial of the matrix $\mbox{\large \tt B}$.
Then it has the following explicit representation$:$
\begin{equation}
\label{P_explicit}
\quad
Q_{n}(\lambda)
=
\left\{
\begin{array}{ccl}
(\lambda-2)\displaystyle\sum\limits_{k=1}^{\lfloor \frac{n\!\!}{\,4} \rfloor +1}{(-1)^{k-1}
{\:\mbox{\scriptsize $\displaystyle\frac{n+1}{2}\!-\!k$}\:
\choose
\:\mbox{\scriptsize $k\!-\!1$}\:}
\lambda^{n-2k+2}}
\!\!&\!\!,\!\!& n\!=\!2m\!+\!1;                                              \\[3.0 ex]
\!\!\!\displaystyle\sum\limits_{k=1}^{\lfloor \frac{n+3}{4} \rfloor +2}{\!\!\!\!(-1)^{k-1}\!{\Bigg (}
\!{\:\mbox{\scriptsize $\displaystyle\frac{n}{2}\!-\!k\!+\!2$}\:
\choose
\mbox{\scriptsize $k\!-\!1$}}
\!+\!
{\:\mbox{\scriptsize $\displaystyle\frac{n}{2}\!-\!k\!+\!2$}\:
\choose
\!\!\mbox{\scriptsize $k\!-\!2$}\:} \! \lambda \! {\Bigg )} \lambda^{n-2k+3}}
\!&\!\!,\!\!& n\!=\!2m.\!\!\!\!
\end{array}
\right.
\end{equation}
\end{lemma}

\smallskip

\noindent
{\bf Proof.}
On the basis of Lemma {\rm \ref{lemma_stara_2}} and Lemma {\rm \ref{lemma_4_2}.} there follows the statement.~\stop

\break

\noindent
The number of the higher order differential operations with Gateaux derivative is determined by
corresponding recurrences, which for dimension $n \!=\! 3, 4, 5,\dots, 10\,$, we can get by the
means of \cite{PowerMatrix}$:$

{\footnotesize
\begin{center}
\begin{tabular}{|c|c|}
\hline $\!\!${\small \rm Dimension}:$\!\!$
&     $\!\!${\small \rm Recurrence$\!$ for$\!$ the$\!$ num.$\!$ of$\!$ the $\!k^{\mbox{\scriptsize \rm th}}\!$-order$\!$ diff.$\!$ operations$\!$
                        with$\!$ Gateaux$\!$ derivative:}$\!\!$     \\ \hline
$n = \;$ 3 & $\mbox{\normalsize \tt g}(k)
             = 2 \mbox{\normalsize \tt g}(k-1)$                                       \\ \hline
$n = \;$ 4 & $\mbox{\normalsize \tt g}(k)
             = \mbox{\normalsize \tt g}(k-1)
             + 2 \mbox{\normalsize \tt g}(k-2)
             - \mbox{\normalsize \tt g}(k-3)$                                         \\ \hline
$n = \;$ 5 & $\mbox{\normalsize \tt g}(k)
             = 2 \mbox{\normalsize \tt g}(k-1)
             + \mbox{\normalsize \tt g}(k-2)
             - 2 \mbox{\normalsize \tt g}(k-3)$                                       \\ \hline
$n = \;$ 6 & $\mbox{\normalsize \tt g}(k)
             = \mbox{\normalsize \tt g}(k-1)
             + 3 \mbox{\normalsize \tt g}(k-2)
             - 2 \mbox{\normalsize \tt g}(k-3)
             - \mbox{\normalsize \tt g}(k-4)$                                         \\ \hline
$n = \;$ 7 & $\mbox{\normalsize \tt g}(k)
             = 2 \mbox{\normalsize \tt g}(k-1)
             + 2 \mbox{\normalsize \tt g}(k-2)
             - 4 \mbox{\normalsize \tt g}(k-3)$                                       \\ \hline
$n = \;$ 8 & $\mbox{\normalsize \tt g}(k)
             = \mbox{\normalsize \tt g}(k-1)
             + 4 \mbox{\normalsize \tt g}(k-2)
             - 3 \mbox{\normalsize \tt g}(k-3)
             - 3 \mbox{\normalsize \tt g}(k-4)
             + \mbox{\normalsize \tt g}(k-5)$                                         \\ \hline
$n = \;$ 9 & $\mbox{\normalsize \tt g}(k)
             = 2 \mbox{\normalsize \tt g}(k-1)
             + 3 \mbox{\normalsize \tt g}(k-2)
             - 6 \mbox{\normalsize \tt g}(k-3)
             - \mbox{\normalsize \tt g}(k-4)
             + 2 \mbox{\normalsize \tt g}(k-5)$                                       \\ \hline
$n =   $10 & $\mbox{\normalsize \tt g}(k)
             = \mbox{\normalsize \tt g}(k-1)
             + 5 \mbox{\normalsize \tt g}(k-2)
             - 4 \mbox{\normalsize \tt g}(k-3)
             - 6 \mbox{\normalsize \tt g}(k-4)
             + 3 \mbox{\normalsize \tt g}(k-5)
             + \mbox{\normalsize \tt g}(k-6)$                                        \\ \hline
\end{tabular}
\end{center}}

\smallskip
\noindent For considered dimensions $n=3,4,5, \dots ,10$, the
values of the function $\mbox{\large \tt g}(k)$, for small values
of the argument $k$, are given in the database of integer sequences
\cite{Sloane07} as sequences \seqnum{A000079} $(n=3)$,
\seqnum{A090990} $(n=4)$, \seqnum{A007283} $(n=5)$, \seqnum{A090992}
$(n=6)$, \seqnum{A000079} $(n=7)$, \seqnum{A090994} $(n=8)$,
\seqnum{A020714} $(n=9)$, \seqnum{A129638} $(n=10)$, respectively.

\bibstyle{plain}

\bigskip
\hrule
\bigskip

\noindent 2000 {\it Mathematics Subject Classification}: 05C30, 26B12, 58C20. \\

\noindent \emph{Keywords: the compositions of the differential operations, enumeration of graphs and maps,
Gateaux directional derivative}

\bigskip
\hrule
\bigskip

\noindent (Concerned with sequence
\seqnum{A000079},
\seqnum{A007283},
\seqnum{A020701},
\seqnum{A020714},
\seqnum{A090989},
\seqnum{A090990},
\seqnum{A090991},
\seqnum{A090992},
\seqnum{A090993},
\seqnum{A090994},
\seqnum{A090995},
\seqnum{A129638})

\vspace*{+.1in} \noindent Received June 5, 2007.
% revised version received  ... .
% Published in {\it Journal of Integer Sequences} ... .
%
% \bigskip
% \hrule
% \bigskip
%
% \noindent
% Return to
% \htmladdnormallink{Journal of Integer Sequences home page}{http://www.math.uwaterloo.ca/JIS/}.

\end{document}